\documentclass[final,3p,times]{elsarticle}

\usepackage{graphicx,amsmath,amssymb,txfonts}
\usepackage{algorithm}
\usepackage{algorithmic}
\usepackage{subfigure,color}
\usepackage[colorlinks]{hyperref}
\usepackage{hypernat}
\usepackage{multirow}
\numberwithin{equation}{section}

\biboptions{sort&compress}

\makeatletter

\newcommand{\Rmnum}[1]{\expandafter\@slowromancap\romannumeral #1@}
\makeatother

\newtheorem{theorem}{\textbf{Theorem}}

\newtheorem{lemma}{\textbf{Lemma}}

\newtheorem{remark}{\textbf{Remark}}
\newtheorem{example}{\textbf{Example}}

\begin{document}

\begin{frontmatter}



\title{A  three-level linearized difference scheme for the coupled nonlinear fractional Ginzburg-Landau equation}


\author[Tongji]{Dongdong He}

\author[CSU]{Kejia Pan\corref{cor}}

\cortext[cor]{Corresponding author.\\ E-mail address: hedongdong@cuhk.edu.cn(D.D. He), pankejia@hotmail.com(K.J. Pan).}
\address[Tongji]{School of Science and Engineering,
The Chinese University of Hong Kong, Shenzhen, Guangdong 518172, China}
\address[CSU]{School of Mathematics and Statistics, Central South University, Changsha, Hunan 410083, China}

\begin{abstract}
In this paper, the coupled fractional Ginzburg-Landau  equations  are first time investigated numerically.  A   linearized implicit finite difference scheme is proposed.   The scheme involves three time levels, is unconditionally stable and second-order accurate in both time and space variables. The unique solvability, the unconditional stability and optimal pointwise error estimates  are obtained by using the energy method and mathematical induction. Moreover, the proposed second-order method can be easily extended into the fourth-order method by using an average finite difference operator for spatial fractional derivatives and Richardson extrapolation for time variable.  Finally, numerical results are presented to confirm the theoretical results.
\end{abstract}

\begin{keyword}
 Ginzburg-Landau  equation, fractional Laplacian, pointwise error estimate, unconditional stability, convergence
\end{keyword}

\end{frontmatter}
 AMS subject classifications: 65M06

\section{Introduction}

The Ginzburg-Landau equation (GLE) has been used to model a wide variety of physical systems~\cite{Aranson}.
The existence, regularity,
and the long-time behavior of the exact solution for the GLE are given in~\cite{DuanJ1,DuanJ2,Doering,Gao}. For
the global well-posedness and the existence of the global attractor of the GLE, readers can refer to Refs.~\cite{GuoB,Huo}. Sakaguchi and Malomed~\cite{Sakaguchi} proposed the coupled Ginzburg-Landau equations (CGLE) to describe the Bose-Einstein condensates and nonlinear optical waveguides and cavities. And the linearly coupled complex cubic-quintic Ginzburg-Landau equations~\cite{Porsezian} and the coupled Ginzburg-Landau equations with higher-order nonlinearities~\cite{Zakeri} are also discussed.

The fractional generalization of the GLE was first suggested by Tarasov and Zaslavsky~\cite{Tarasov1,Tarasov2} for fractal media. Since then, the fractional Ginzburg-Landau equation (FGLE) has been exploited to describe various physical phenomena, such as the media with fractal mass dimension~\cite{Tarasov1}, a class of critical phenomena when the organization of the system near the phase transition point is influenced by a competing nonlocal ordering~\cite{Milovanov} and a network of diffusively Hindmarsh-Rose neurons with long-range synaptic coupling~\cite{Mvogo}.  Recently, the coupled fractional Ginzburg-Landau equation (CFGLE) with stocastic noise was discussed by~\cite{Shu}.

The dissipative mechanism of the FGLE and CFGLE are not characterized by the classical Laplacian but by the fractional power of the Laplacian, which brings some essential difficulties in theoretical analysis. Many authors have worked on the FGLE and CFGLE from the theoretical aspects~\cite{Shu,Tarasov,Pu,Guo,Lu,Millot}.
For example, the solution of FGLE is derived by using psi-series with fractional powers~\cite{Tarasov}.
 Pu and Guo~\cite{Pu} investigated the global well-posedness and dynamics for the nonlinear FGLE. Millot and Sire~\cite{Millot} considered the asymptotic analysis of the FGLE in a bounded domain. Shu {\it et al.}~\cite{Shu} studied the random attractors for the CFGLE with stochastic noise.

From the numerical point of view, there are quite a lot of numerical studies for the classical GLE, see~\cite{Lord,Xu,Wang,ZhangL,ZhangYN,Hao} and references therein. But there are not too much numerical studies for the FGLE.
Mvogo et al.~\cite{Mvogo} proposed a semi-implicit Riesz fractional finite difference scheme, which is only first-order accurate in time and second-order accurate in space. Wang and Huang~\cite{WangP1} proposed an implicit midpoint difference scheme for the FGLE, which is shown to be second-order accurate in the $L^2$-norm for  both time and space variables.
However, the method is a nonlinear scheme, a fixed point iteration is needed at each time step, which is generally computational expensive.
Hao and Sun~\cite{Hao2016} proposed a three-level linearized method for FGLE with second-order accuracy in time and fourth-order accuracy in space. However, this method is not unconditionally stable. The authors~\cite{He2018} recently proposed  an unconditionally stable three-level linearized scheme  for the FGLE with second-order accuracy in both time and space variables, and method in~\cite{He2018} can be easily extended into fourth-order method in  space variable with unconditional stability. Very recently, Wang and Huang~\cite{Wang2018} proposed another three-level linearized difference method for FGLE.  To our best knowledge, although there are some theoretical studies for the CGLE and CFGLE, there is no numerical investigation.

In this paper, we consider the following CFGLE
\begin{align}
u_t+(\upsilon_1+i\eta_1)(-\Delta)^{\frac{\alpha}{2}}u+\left((\kappa_1+i\zeta_1)|u|^2+(\delta_1+i\beta_1)|v|^2\right)u-\gamma_1 u&=0, \quad x\in \mathbb{R}, \quad 0<t<T,\label{IVP1}\\
v_t+(\upsilon_2+i\eta_2)(-\Delta)^{\frac{\alpha}{2}}v+\left((\kappa_2+i\zeta_2)|u|^2+(\delta_2+i\beta_2)|v|^2\right)v-\gamma_2 v&=0, \quad x\in \mathbb{R}, \quad 0<t<T, \label{IVP2}
\end{align}
with initial condition
\begin{align}\label{IVP3}
u(x,0)=u_0(x),\quad v(x,0)=v_0(x), \quad x\in \mathbb{R}, \quad 0<t<T,
\end{align}
where $i=\sqrt{-1}$ is the complex unit, $u(x, t), v(x,t)$ are complex-valued functions of time variable $t$ and space variable $x$, $\upsilon_1>0,\upsilon_2>0,\kappa_1,\kappa_2,\delta_1,\delta_2,\eta_1,\eta_2,\zeta_1,\zeta_2,\gamma_1,\gamma_2$ are given real constants, $u_0(x), v_0(x)$ are   complex-valued functions satisfying certain regularity, and $1<\alpha\leq 2$. The fractional Laplacian can be regarded as a pseudo-differential operator with the symbol $-|\xi|^{\alpha}$:
$$-(-\Delta)^{\frac{\alpha}{2}}u(x,t)=\mathcal{F}(|\xi|^{\alpha}\hat{u}(\xi,t)),$$
where $\mathcal{F}$ denotes the Fourier transform. It is indeed equivalent to the following Riesz fractional derivative~\cite{Yang,Zhuang,Wang2013,Podlubny}
\begin{align*}
-(-\Delta)^{\frac{\alpha}{2}}u=\frac{\partial^{\alpha}u(x,t)}{\partial|x|^{\alpha}}:=-\frac{1}{2\cos(\pi\frac{\alpha}{2})}[_{-\infty}D^{\alpha}_{x}u(x,t)+_{x}D^{\alpha}_{+\infty}u(x,t)],
\end{align*}
where $_{-\infty}D^{\alpha}_{x}$ and $_{x}D^{\alpha}_{+\infty}$ are left- and right-side Riemann-Liouville fractional derivatives, respectively. {Some other definitions of fractional derivatives can be found in~\cite{Podlubny,Xu2013}.
Obviously, when $\alpha=2$, the fractional Laplacian reduces to the classical Laplacian and Eqs. (\ref{IVP1})-(\ref{IVP2}) reduce to the coupled nonlinear GLE. Furthermore, if $\upsilon_1 = \kappa_1=\delta_1=\gamma_1=\upsilon_2 = \kappa_2=\delta_2=\gamma_2=0$,
Eqs. (\ref{IVP1})-(\ref{IVP2}) reduce to the coupled nonlinear Schr${\ddot{\rm o}}$dinger equation~\cite{Wang2013}.}
The Eqs. (\ref{IVP1})-(\ref{IVP2}) can also be viewed as the direct generalization of the FGLE discussed in~\cite{He2018}.

The objective of this paper is to develop an unconditionally stable linearized scheme with optimal pointwise error estimates for the  initial boundary value problem (\ref{IVP1})-(\ref{IVP2}).
 The method, which uses three time levels, is shown to be second-order convergent in the $L^{\infty}$-norm for both  time and space variables.  And the method  is also shown to be almost unconditionally stable (the time step is not related to the spatial meshsize). It should be noted that, combined with  {high-order} difference schemes for the spatial approximation, our method can be easily extended to spatial fourth-order accuracy with almost unconditional stability.

The rest of this paper is organized as follows. Section~\ref{section2} reviews the central difference method for the spatial fractional derivatives and the fractional Sobolev space. Section~\ref{section3} gives the linearized implicit finite difference method. Section~\ref{section4} provides the theoretical analysis for the proposed scheme, which includes the unique solvability, convergence and stability. Section~\ref{section5} presents the  numerical results which confirm the theoretical results. And the conclusion is given in the final section.

\section{Preliminaries}\label{section2}
\subsection{Spatial discretization}

For approximation of the Riesz fractional derivative, considerable efforts have been made. For example, Meerschaert and Tadjeran~\cite{Meerschaert1,Meerschaert2,Meerschaert3} proposed the shifted Gr$\ddot{\rm u}$nwald formula. Yang et al.~\cite{Yang} developed three numerical methods, namely, the L1/L2 approximation method, the standard/shifted Gr$\ddot{\rm u}$nwald formula, and the matrix transform method~\cite{Ilic,Yang2}. Zhou et al.~\cite{Tian,Zhou} proposed a weighted and shifted Gr$\ddot{\rm u}$nwald difference scheme. A finite element method was presented in~\cite{Zhang}. Recently, the fractional centered difference was defined by Ortigueira~\cite{Ortigueira}. $\c{\rm C}$elik and Duman~\cite{Celik} analyzed the error of this approximation and applied it to fractional diffusion equations. Shen et al.~\cite{Shen} employed two fractional centered differences and then proposed a weighted difference scheme for the Riesz space fractional advection-dispersion equation. In this paper, we will adopt this fractional centered difference discretization for the Riesz fractional derivative.

For $\alpha>-1$, the fractional centered difference is defined by~\cite{Ortigueira}
\begin{equation}\label{abc}
\Delta^{\alpha}_h f(x)=\frac{1}{h^{\alpha}}\sum^{\infty}_{k=-\infty}c^{\alpha}_{k}f(x-kh).
\end{equation}
where
\begin{equation}
  c^{\alpha}_k=\frac{(-1)^{k}\Gamma(\alpha+1)}{\Gamma(\frac{\alpha}{2}-k+1)\Gamma(\frac{\alpha}{2}+k+1)}=\Big(1-\frac{\alpha+1}{\frac{\alpha}{2}+k}\Big)c_{k-1}^\alpha, \quad \textrm{for } k\in \mathbb{Z}.
\end{equation}



 \begin{lemma}\label{Lemma1}
  (see~\cite{Celik}) The coefficients $c^{\alpha}_k$  have the following properties for $-1<\alpha\leq 2$
\begin{equation}
c^{\alpha}_0>0,\quad c^{\alpha}_k=c^{\alpha}_{-k}\leq 0,\quad \textrm{for}\ k=\pm1, \pm 2, \cdots.
\end{equation}
 \end{lemma}

 \begin{lemma}\label{Lemma1_1}
  (see~\cite{Celik}) Let $f(x)\in C^{5}(\mathbb{R})\cap L^1(\mathbb{R})$ and all spatial derivatives of $f(x)$ up to order five belong to $L^1(\mathbb{R})$. Then
\begin{equation}
(-\Delta)^{\frac{\alpha}{2}}f(x) = \Delta^{\alpha}_h f(x) + O(h^2),
\end{equation}
for $1<\alpha \leq 2$.
 \end{lemma}

\subsection{Fractional Sobolev norm}

Denote infinite grid by $Z_h$ with grid points $x_j=jh$ ($j\in Z$). For any grid functions {$U=\{U_j\}, V=\{V_j\}$} on $Z_h$, the discrete inner product and the norms are defined as
\begin{align}\label{notation2}
 {(U,V)=h\sum_{j\in Z}U_j\bar{V}_j,\quad \parallel U\parallel=\sqrt{(U,U)_h}, \quad \parallel U\parallel_{\infty}=\max_{j\in Z}{|U_j|}.}
\end{align}

Set {$L^2_h=\left\{V|V\in Z_h,\parallel V\parallel<+\infty\right\}$}, then for a given constant $\sigma \in[0,1]$, the fractional Sobolev norm $\parallel V\parallel_{H^{\sigma}}$  and seminorm $|V|_{H^{\sigma}}$ can be defined as follows~\cite{WangP1}
\begin{align}\label{Sobolev}
\parallel V\parallel^2_{H^{\sigma}}=\int^{\pi/h}_{-\pi/h}(1+|k|^{2\sigma})|\hat{V}(k)|^2dk,\quad |V|^2_{H^{\sigma}}=\int^{\pi/h}_{-\pi/h}|k|^{2\sigma}|\hat{V}(k)|^2dk,
\end{align}
where the semi-discrete Fourier transform $\hat{V}(k)$ is defined as
\begin{align}\label{Fourier}
\hat{V}(k)=\frac{h}{\sqrt{2\pi}}\sum_{j\in Z}V_je^{-ikx_j}.
\end{align}
Moreover, we also have the inversion formula
\begin{equation}
  V_j = \frac{1}{\sqrt{2\pi}}\int_{-\pi/h}^{\pi/h} \hat{V}(k)e^{ikx_j}dk.
\end{equation}
Obviously, $\parallel V\parallel^2_{H^{\sigma}}=\parallel V\parallel^2+|V|^2_{H^{\sigma}}$. Set {$H_h^\sigma:=\left\{V|V\in Z_h, ||V||_{H^{\sigma}}<+\infty\right\}$.}
%

 \begin{lemma}\label{Lemma2}
 (Discrete Sobolev Inequality~\cite{Huang2016b}) For any $1<\alpha\leq 2$ and $v\in H_h^{\frac{\alpha}{2}}$, it holds that
 \begin{align}
\parallel v\parallel_{\infty}\leq C_{\alpha}\parallel v\parallel_{H^{\frac{\alpha}{2}}},
\end{align}
where $C_{\alpha}>0$ is a constant which depends on $\alpha$ but independent of $h$.
\end{lemma}
%

 \begin{lemma}\label{Lemma3}
 (Fractional seminorm equivalence~\cite{Huang2016b}) For every $1<\alpha\leq 2$, we have
 \begin{align}
\left(\frac{2}{\pi}\right)^{\alpha}|U|^2_{H^{\frac{\alpha}{2}}}\leq (\Delta^{\alpha}_{h}U,U)\leq |U|^2_{H^{\frac{\alpha}{2}}}, \quad \textrm{for\ any}\ U\in {H_h^\frac{\alpha}{2}}.
\end{align}
In addition, the following inequality is also valid,
 \begin{align}\label{Cauchyinquality}
\left(\frac{2}{\pi}\right)^{\alpha}|U|_{H^{\frac{\alpha}{2}}}|V|_{H^{\frac{\alpha}{2}}}\leq \sum_{j\in Z}|\Delta^{\alpha}_{h}U_j\bar{V}_j|\leq |U|_{H^{\frac{\alpha}{2}}} |V|_{H^{\frac{\alpha}{2}}}, \quad \textrm{for\ any}\ {U, V\in H_h^\frac{\alpha}{2}}.
\end{align}
\end{lemma}

\section{A three-level linearized implicit difference scheme}\label{section3}




Under certain conditions, the solution of problem (\ref{IVP1})-(\ref{IVP3}) will be small when $|x|\rightarrow\infty$. Thus, in practical numerical computation, we can truncate the original problem on a bounded interval and take the following homogeneous Dirichlet boundary conditons:
\begin{align}
u(a,t)=0,\quad u(b,t)=0,  \quad v(a,t)=0,\quad v(b,t)=0, \quad t\in [0,T],
\end{align}
where $a$ and $b$ are usually chosen sufficiently large negative and positive numbers such that the truncation error is negligible.
{In the following, we will not make an explicit distinction between the solution $u(x,t)$ defined in $\mathbb{R}$ and the solution $u(x,t)$ defined in $[a,b]$.}

Now we give the numerical discretization  of the homogeneous boundary value problem (\ref{IVP1})-(\ref{IVP2}) in the finite domain $\Omega=[a,b]$. The solution domain is defined as $\{(x,t) | a\leq x\leq b,\ 0\leq t\leq T\}$, which is covered by a uniform grid $\{(x_j,t_n)|x_j=a+jh,\ t_n=n\tau,\ j=0, \cdots, M,\ n=0, \cdots, N\}$, with spacing $h=\frac{b-a}{M},\ \tau=\frac{T}{N}$, where $M,N$ are two positive integers.

Let $U^n_j, V^n_j$ denote the numerical approximation to $u(x_j,t_n)$ and $v(x_j,t_n)$, respectively.  For any grid function $U^n=(U_0^n,U_1^n,\cdots,U_M^n)$, $V^n=(V_0^n,V_1^n,\cdots,V_M^n)$, we define
 \begin{align}
 U^{\bar{n}}_{j}=\frac{U^{n+1}_{j}+U^{n-1}_{j}}{2}, \quad \delta_{t}U^{n}_j=\frac{U^{n+1}_j-U^{n-1}_j}{2\tau}, \quad  V^{\bar{n}}_{j}=\frac{V^{n+1}_{j}+V^{n-1}_{j}}{2}, \quad \delta_{t}V^{n}_j=\frac{V^{n+1}_j-V^{n-1}_j}{2\tau}.
 \end{align}
 Denote$$Z^{0}_h=\left\{V|V=\{V_j\}, j=0, 1, \cdots M, V_{0}=V_{M}=0\right\}.$$
For any $U, V \in Z^{0}_h$, {the discrete} inner product and norms  are defined as follows:
\begin{align}\label{notation1}
(U,V)=h\sum^{M-1}_{j=1}U_j\bar{V}_j,\quad \parallel U\parallel=\sqrt{(U,U)},  \quad \parallel U\parallel_{\infty}=\max_{1\leq j \leq M-1}{|U_j|},
\end{align}
where $\bar{V}$ is the conjugate of $V$.


We note that for any $U,V\in Z^{0}_h$, we can extend the $U, V$ into an infinity grid function by assigning $U_j=0, V_j=0$ for $j\neq 0, 1, \cdots,M$, then the  norms and inner product defined in (\ref{notation2}) can also be defined for grid functions in $Z^{0}_h$. Moreover, Lemma~\ref{Lemma2} and Lemma~\ref{Lemma3}  are also valid for any $U\in Z^{0}_h$. Here and after we don't make a distinction between the inner products and norms defined in (\ref{notation2}) and (\ref{notation1}).

With the assumption of homogenous boundary condition (\ref{IVP3}),
for any $U\in Z^{0}_h$, we have
\begin{equation}
\Delta^{\alpha}_hU^{}_{j}=\frac{1}{h^{\alpha}}\sum^{j-M+1}_{k=j-1}c^{\alpha}_{k}U^{}_{j-k}=\frac{1}{h^{\alpha}}\sum^{M-1}_{k=1}c^{\alpha}_{j-k}U^{}_{k}.
\end{equation}
By using three time levels, our linearized implicit method with for the homogeneous boundary value problem (\ref{IVP1})-(\ref{IVP2}) in the finite domain $\Omega=[a,b]$ is as follows,
\begin{align}
\delta_{t}U^{n}_j+(\upsilon_1+i\eta_1)\Delta^{\alpha}_hU^{\bar{n}}_{j}+\left((\kappa_1+i\zeta_1)|U^{n}_{j}|^2+(\delta_1+i\beta_1)|V^{n}_{j}|^2\right)U^{\bar{n}}_{j}-\gamma_1 U^{\bar{n}}_{j}&=0,\quad 0<j<M,\ 1<n<N,\label{diff1}\\
\delta_{t}V^{n}_j+(\upsilon_2+i\eta_2)\Delta^{\alpha}_hV^{\bar{n}}_{j}+\left((\kappa_2+i\zeta_2)|U^{n}_{j}|^2+(\delta_2+i\beta_2)|V^{n}_{j}|^2\right)V^{\bar{n}}_{j}-\gamma_2 V^{\bar{n}}_{j}&=0,\quad 0<j<M,\ 1<n<N.\label{diff1l}
\end{align}
with initial condition
\begin{equation}\label{diff2}
U^{0}_j=(u_0)_j=u_0(x_j),\quad V^{0}_j=(v_0)_j=v_0(x_j), \quad 0<j<M,
\end{equation}
and boundary conditions
\begin{equation}\label{diff3}
U^{n}_0=U^{n}_M=0,\quad V^{n}_0=V^{n}_M=0,\quad 0\leq n\leq N.
\end{equation}

Since the difference scheme (\ref{diff1}) involves three time levels, the first step values $U^1_j, V^1_j$ are required to begin stepping the numerical solution forward in time.
From Taylor expansion, one has
\begin{align}\label{firststep}
u^1=&u_0+\tau u_t(x,0)+\frac{\tau^2}{2}u_{tt}(x,\bar{\tau})=u_0-\tau\left((\upsilon_1+i\eta_1)(-\Delta)^{\frac{\alpha}{2}}u_0+\left((\kappa_1+i\zeta_1)|u_0|^2+(\delta_1+i\beta_1)|v_0|^2\right)u_0-\gamma_1 u_0\right)+\frac{\tau^2}{2}u_{tt}(x,\bar{\tau}_1),\\
v^1=&v_0+\tau v_t(x,0)+\frac{\tau^2}{2}v_{tt}(x,\bar{\tau})=v_0-\tau\left((\upsilon_2+i\eta_2)(-\Delta)^{\frac{\alpha}{2}}v_0+\left((\kappa_2+i\zeta_2)|u_0|^2+(\delta_2+i\beta_2)|v_0|^2\right)v_0-\gamma_2 v_0\right)+\frac{\tau^2}{2}v_{tt}(x,\bar{\tau}_2),
\end{align}
{for $\bar{\tau}_1, \bar{\tau}_2$ satisfying $0<\bar{\tau}_1, \bar{\tau}_2<\tau$ and}
\begin{align*}
   u_t &= -(\upsilon_1+i\eta_1)(-\Delta)^{\frac{\alpha}{2}}u-\left((\kappa_1+i\zeta_1)|u|^2+(\delta_1+i\beta_1)|v|^2\right)u+\gamma_1 u,\\ v_t &= -(\upsilon_2+i\eta_2)(-\Delta)^{\frac{\alpha}{2}}v-\left((\kappa_2+i\zeta_2)|u|^2+(\delta_2+i\beta_2)|v|^2\right)v+\gamma_2 v,
 \end{align*}
are used in the above two equalities.

In addition, from Lemma \ref{Lemma1_1} we know
\begin{align}\label{secondstep}
(-\Delta)^{\frac{\alpha}{2}}(u_0)_j&=\frac{1}{h^{\alpha}}\sum^{M-1}_{k=1}c^{\alpha}_{j-k}u^{0}_{k}+O(h^2)=\Delta^{\alpha}_hu^{0}_{j}+O(h^2)\\
(-\Delta)^{\frac{\alpha}{2}}(v_0)_j&=\frac{1}{h^{\alpha}}\sum^{M-1}_{k=1}c^{\alpha}_{j-k}v^{0}_{k}+O(h^2)=\Delta^{\alpha}_hv^{0}_{j}+O(h^2).
\end{align}
In the numerical simulation, $U^1, V^1$ are obtained from the following scheme
\begin{align}
U^1_j&=U^0_j-\tau\left((\upsilon_1+i\eta_1)\Delta^{\alpha}_hU^0_j+\left((\kappa_1+i\zeta_1)|U^0_j|^2+(\delta_1+i\beta_1)|V^0_j|^2\right)U^0_j-\gamma_1 U^0_j\right),\quad 0<j<M,\label{diffinitialaa}\\
V^1_j&=V^0_j-\tau\left((\upsilon_2+i\eta_2)\Delta^{\alpha}_hV^0_j+\left((\kappa_2+i\zeta_2)|U^0_j|^2+(\delta_2+i\beta_2)|V^0_j|^2\right)V^0_j-\gamma_2 V^0_j\right),\quad 0<j<M.\label{diffinitialaabb}
\end{align}

\section{Theoretical analysis}\label{section4}

%

The following four lemmas are essential for the analysis of the numerical solution.

 \begin{lemma}\label{Lemma4}
 (see~\cite{WangP2}) For any two grid functions $U,V\in Z^{0}_h$, there exists a linear operator $\Lambda^{\alpha}$ such that
 \begin{align}
(\Delta^{\alpha}_{h}U,V)=(\Lambda^{\alpha}U,\Lambda^{\alpha}V).
\end{align}
\end{lemma}

By using the above lemma, the following lemma is easy to verify.
 \begin{lemma}\label{Lemma5}
 For any grid functions $U^n\in Z^{0}_h$, we have
 \begin{align}
&Im(\Delta^{\alpha}_{h}U^n,U^n)=0,\\
&Re(\Delta^{\alpha}_{h}U^{\bar{n}},\delta_{t} U^n)=\frac{1}{4\tau}(\parallel\Lambda^{\alpha} U^{n+1}\parallel^2-\parallel\Lambda^{\alpha} U^{n-1}\parallel^2).
\end{align}
\end{lemma}

{
 \begin{lemma}\label{Lemma6}
 (Discrete Gronwall's inequality~\cite{Holte,Hu2015}). Let $\{u_k\}$ and $\{w_k\}$ be nonnegative sequences and $\alpha$ a nonnegative constant satisfying
 \begin{equation}
   u_n\leq \alpha + \sum_{0\leq k<n}w_k u_k \quad \textrm{for } n\geq 0.
 \end{equation}
Then for all $n$ it holds
\begin{equation}
   u_n\leq \alpha \exp\left(\sum_{0\leq k<n} w_k\right).
\end{equation}
\end{lemma}}

\subsection{$L^{\infty}$ convergence and stability}

Let $u(x,t), v(x,t)$ be the exact solution of problem (\ref{IVP1})-(\ref{IVP3}),  $U^n_j, V^n_j$ be the solution of the numerical schemes (\ref{diff1})-(\ref{diff3}). Let  $u^n_j=u(x_j,t_n), v^n_j=v(x_j,t_n)$, the error functions $$e^n_j=u^n_j-U^n_j,\quad \xi^n_j=v^n_j-V^n_j,\quad j=1,2,\cdots, M,\;\; n=1,2,\cdots,N.$$

Define the truncation errors of the scheme (\ref{diff1})-(\ref{diff1l}) as follows:
\begin{align}
r^n_j & =\delta_{t}u^{n}_j+(\upsilon_1+i\eta_1)\Delta^{\alpha}_hu^{\bar{n}}_{j}+\left((\kappa_1+i\zeta_1)|u^{n}_{j}|^2+(\delta_1+i\beta_1)|v^{n}_{j}|^2\right)u^{\bar{n}}_{j}-\gamma_1 u^{\bar{n}}_{j},\quad 1\leq j\leq M-1,\quad  1\leq n \leq N-1,\label{trun1}\\
s^n_j &=\delta_{t}v^{n}_j+(\upsilon_2+i\eta_2)\Delta^{\alpha}_hv^{\bar{n}}_{j}+\left((\kappa_2+i\zeta_2)|u^{n}_{j}|^2+(\delta_2+i\beta_2)|v^{n}_{j}|^2\right)v^{\bar{n}}_{j}-\gamma_2 v^{\bar{n}}_{j},\quad 1\leq j\leq M-1,\quad  1\leq n \leq N-1.\label{trun2}
\end{align}

Subtracting (\ref{diff1}) from (\ref{trun1}) and subtracting (\ref{diff1l}) from (\ref{trun2}) yield that
\begin{align}
r^n_j=&\delta_{t}e^n_j+(\upsilon_1+i\eta_1)\Delta^{\alpha}_he^{\bar{n}}_{j}+P^n_j-\gamma_1 e^{\bar{n}}_{j},\label{truncation}\\
s^n_j=&\delta_{t}\xi^n_j+(\upsilon_2+i\eta_2)\Delta^{\alpha}_h\xi^{\bar{n}}_{j}+Q^n_j-\gamma_2 \xi^{\bar{n}}_{j},\label{truncation0}
\end{align}
where $$P^n_j=\left((\kappa_1+i\zeta_1)|u^n_j|^2+(\delta_1+i\beta_1)|v^n_j|^2\right)u^{\bar{n}}_{j}-\left((\kappa_1+i\zeta_1)|U^n_j|^2+(\delta_1+i\beta_1)|V^n_j|^2\right)U^{\bar{n}}_{j},$$
$$Q^n_j=\left((\kappa_2+i\zeta_2)|u^n_j|^2+(\delta_2+i\beta_2)|v^n_j|^2\right)v^{\bar{n}}_{j}-\left((\kappa_2+i\zeta_2)|U^n_j|^2+(\delta_2+i\beta_2)|V^n_j|^2\right)V^{\bar{n}}_{j},$$
for $1\leq j\leq M-1$ and $1\leq n \leq N-1$.

Using Taylor expansion and Lemma \ref{Lemma1_1}, we can easily obtain the following lemma.
 \begin{lemma}\label{Lemma8}
Suppose that the solution of problem (\ref{IVP1})-(\ref{IVP3}) is sufficiently smooth. Then it is holds that
\begin{align}
&|r^n_j|\leq C_R(\tau^2+h^2),\quad |s^n_j|\leq C_R(\tau^2+h^2), \quad 1\leq j\leq M-1,\quad  1\leq n \leq N-1,
\end{align}
where $C_R$ is a positive constant independent of $\tau$ and $h$.
\end{lemma}

Following a similar proof of Lemma 9 in \cite{He2018}, we can obtain the lemma below.
 \begin{lemma}\label{Lemma11}
 Suppose that the solution of problem (\ref{IVP1})-(\ref{IVP3}) is sufficiently smooth. Then one has
 \begin{equation}\label{e1}
   \begin{aligned}
|e^1_j|&\leq C_e(\tau^2+\tau h^2),\quad  |\Delta^{\alpha}_h e^1_j|\leq C_e(\tau^2+\tau h^{2+\alpha}),\\
 |\xi^1_j|&\leq C_e(\tau^2+\tau h^2),\quad  |\Delta^{\alpha}_h \xi^1_j|\leq C_e(\tau^2+\tau h^{2+\alpha}),
\end{aligned}
 \end{equation}
where $C_e$ is a positive constant independent of $\tau$ and $h$.
\end{lemma}

 \begin{theorem}\label{theorem3}
 Suppose that the solution of problem (\ref{IVP1})-(\ref{IVP3}) is smooth enough, then there exist two small positive constants $\tau_0$ and $h_0$,  such that, when $\tau<\tau_0$ and $h<h_0$, the numerical solution $(U^n, V^n)$ of the difference scheme (\ref{diff1})-(\ref{diff3}) converges to the exact solution $(u^n, v^n)$ in the sense of $L^{\infty}$-norm with the optimal convergence order $O(\tau^2+h^2)$, i.e.,
 \begin{align}
\parallel u^n-U^n\parallel_{\infty}&\leq C_0(\tau^2+h^2), \quad \parallel v^n-V^n\parallel_{\infty}\leq C_0(\tau^2+h^2), \quad  1\leq n \leq N, \label{estimate1}
\end{align}
where $C_0$ is a positive constant independent of $\tau$ and $h$.
 \end{theorem}

\noindent {\bf Proof.}
we use mathematical induction to prove (\ref{estimate1}). It follows from (\ref{e1}) that
the error estimate (\ref{estimate1}) holds for $n=1$.
Assume that  (\ref{estimate1}) is valid for $m\leq n$, we want to show that (\ref{estimate1}) is also valid for $n+1$.

By the assumption, one has
\begin{equation}\label{estimate47}
  \begin{aligned}
\parallel U^{m}\parallel_{\infty}\leq \parallel u^{m}\parallel_{\infty}+\parallel e^{m}\parallel_{\infty}\leq C_m + C_0(\tau^2+h^2) \leq  C_m+1, \quad  1\leq m \leq n, \\
\parallel V^{m}\parallel_{\infty}\leq \parallel v^{m}\parallel_{\infty}+\parallel \xi^{m}\parallel_{\infty}\leq C_m + C_0(\tau^2+h^2) \leq  C_m+1, \quad  1\leq m \leq n,
\end{aligned}
\end{equation}
for $\tau<\tau_1$ and $h<h_1$, where $\tau_1, h_1$  satisfy that  $\tau_1^2+h_1^2<\frac{1}{C_0}$. Here,
 $$C_m=\max\left\{\max_{a\leq x \leq b, 0\leq t \leq T }|u(x,t)|, \max_{a\leq x \leq b, 0\leq t \leq T }|v(x,t)|\right\}.$$

Now computing the discrete inner product of (\ref{truncation}) with $e^{\bar{n}}$ and taking the real part of the resulting equation, we have
  \begin{align}\label{estimate48}
\frac{\parallel e^{n+1}\parallel^2-\parallel e^{n-1}\parallel^2}{4\tau}+\upsilon_1\parallel\Lambda^{\alpha}e^{\bar{n}}\parallel^2+Re\left[(P^n,e^{\bar{n}})\right]-\gamma_1\parallel e^{\bar{n}}\parallel^2=Re\left[(r^n,e^{\bar{n}})\right],
\end{align}
where
\begin{align}
P^n_j&=(\kappa_1+i\zeta_1)|u^{n}_{j}|^2u^{\bar{n}}_{j}-(\kappa_1+i\zeta_1)|U^{n}_{j}|^2U^{\bar{n}}_{j}+(\delta_1+i\beta_1)|v^{n}_{j}|^2u^{\bar{n}}_{j}-(\delta_1+i\beta_1)|V^{n}_{j}|^2U^{\bar{n}}_{j}\nonumber\\
     &=(\kappa_1+i\zeta_1)\left((|u^{n}_{j}|^2-|U^{n}_{j}|^2)u^{\bar{n}}_{j}+|U^{n}_{j}|^2e^{\bar{n}}_{j}\right)+(\delta_1+i\beta_1)\left((|v^{n}_{j}|^2-|V^{n}_{j}|^2)u^{\bar{n}}_{j}+|V^{n}_{j}|^2e^{\bar{n}}_{j}\right)\nonumber\\
     &=(\kappa_1+i\zeta_1)\left((|u^{n}_{j}|-|U^{n}_{j}|)(|u^{n}_{j}|+|U^{n}_{j}|)u^{\bar{n}}_{j}+|U^{n}_{j}|^2e^{\bar{n}}_{j}\right)+(\delta_1+i\beta_1)\left(|v^{n}_{j}|-|V^{n}_{j}|)(|v^{n}_{j}|+|V^{n}_{j}|)u^{\bar{n}}_{j}+|V^{n}_{j}|^2e^{\bar{n}}_{j}\right).
\end{align}

By using the assumption (\ref{estimate47}), one has
\begin{align}
\parallel P^n\parallel\leq& \sqrt{\kappa^2_1+\xi^2_1}\left((\parallel u^{n}\parallel_{\infty}+\parallel U^{n}\parallel_{\infty})\parallel u^{\bar{n}}\parallel_{\infty}\parallel e^n\parallel+\parallel U^{n}\parallel^2_{\infty}\parallel e^{\bar{n}}\parallel\right)\nonumber\\
&+\sqrt{\delta^2_1+\beta^2_1}\left((\parallel v^{n}\parallel_{\infty}+\parallel V^{n}\parallel_{\infty})\parallel u^{\bar{n}}\parallel_{\infty}\parallel \xi^n\parallel+\parallel V^{n}\parallel^2_{\infty}\parallel e^{\bar{n}}\parallel\right)\nonumber\\
\leq& \sqrt{\kappa^2_1+\xi^2_1}\left((C_m+C_m+1)(C_m+1)\parallel e^{n}\parallel+(C_m+1)^2\parallel e^{\bar{n}}\parallel\right)\nonumber\\
&+\sqrt{\delta^2_1+\beta^2_1}\left((C_m+C_m+1)(C_m+1)\parallel \xi^{n}\parallel+(C_m+1)^2\parallel e^{\bar{n}}\parallel\right)\nonumber\\
\leq& C_1 (\parallel e^{n}\parallel+\parallel \xi^{n}\parallel+\parallel e^{\bar{n}}\parallel), \label{estimate49}
\end{align}
for $\tau<\tau_1$ and $h<h_1$, where
\begin{equation*}
  C_1 = 2(C_m+1)^2\max\left\{\sqrt{\kappa^2_1+\xi^2_1},\sqrt{\delta^2_1+\beta^2_1}\right\}.
\end{equation*}

Moreover,
\begin{align}
\left|Re\left[(P^n,e^{\bar{n}})\right]\right|&\leq\left|(P^n,e^{\bar{n}})\right|\nonumber\\
&\leq\parallel P^n\parallel\ \parallel e^{\bar{n}}\parallel\nonumber\\
&\leq\frac{1}{2}\left(\parallel P^n\parallel^2+\parallel e^{\bar{n}}\parallel^2\right)\nonumber\\
&\leq \frac{3C^2_1+1}{2}\left(\parallel e^{n-1}\parallel^2+\parallel e^n \parallel^2+\parallel e^{n+1} \parallel^2+\parallel \xi^{n}\parallel^2\right). \label{dd}
\end{align}
Similarly, one has
\begin{align}\label{ee}
\left|Re(r^n,e^{\bar{n}})\right|&\leq|(r^n,e^{\bar{n}})|\leq \frac{1}{2}\left(\parallel r^n \parallel^2+\parallel e^{\bar{n}}\parallel^2\right)\leq \frac{1}{2}\left(\parallel r^n \parallel^2+\parallel e^{n+1}\parallel^2+\parallel e^{n-1}\parallel^2\right).
\end{align}

From (\ref{estimate48}), (\ref{dd}) and (\ref{ee}), we have
\begin{align}\label{estimate50}
&\frac{\parallel e^{n+1}\parallel^2-\parallel e^{n-1}\parallel^2}{4\tau}\\
=&-\upsilon_1\parallel\Lambda^{\alpha}e^{\bar{n}}\parallel^2-Re\left[(P^n,e^{\bar{n}})\right]+\gamma_1\parallel e^{\bar{n}}\parallel^2+Re(r^n,e^{\bar{n}})\nonumber \\
\leq &-Re\left[(P^n,e^{\bar{n}})\right]+\gamma_1\parallel e^{\bar{n}}\parallel^2+Re(r^n,e^{\bar{n}})\nonumber\\
\leq &\frac{3C^2_1+1}{2}(\parallel e^{n+1}\parallel^2+\parallel e^{n}\parallel^2+\parallel e^{n-1}\parallel^2+\parallel \xi^{n}\parallel^2)+\frac{|\gamma_1|}{2}(\parallel e^{n+1}\parallel^2+\parallel e^{n-1}\parallel^2)+\frac{1}{2}(\parallel e^{n+1}\parallel^2+\parallel e^{n-1}\parallel^2+\parallel r^n \parallel^2)\nonumber\\
\leq& C_2(\parallel e^{n+1}\parallel^2+\parallel e^{n}\parallel^2+\parallel e^{n-1}\parallel^2+\parallel \xi^{n}\parallel^2)+\frac{1}{2}\parallel r^n \parallel^2,
\end{align}
where $C_2=\frac{3C^2_1+|\gamma_1|+2}{2}$.
Thus,
\begin{align}\label{estimate51}
\parallel e^{n+1}\parallel^2-\parallel e^{n-1}\parallel^2\leq  4C_2\tau(\parallel e^{n+1}\parallel^2+\parallel e^{n}\parallel^2+\parallel e^{n-1}\parallel^2+\parallel \xi^{n}\parallel^2)+2\tau\parallel r^n\parallel^2.
\end{align}

%
%
Similarly,  computing the discrete inner product of (\ref{truncation0}) with $\xi^{\bar{n}}$ and analyzing the resulting equation, one  can obtain
\begin{align}\label{estimate510}
\parallel \xi^{n+1}\parallel^2-\parallel \xi^{n-1}\parallel^2\leq  4C'_2\tau(\parallel \xi^{n+1}\parallel^2+\parallel \xi^{n}\parallel^2+\parallel \xi^{n-1}\parallel^2+\parallel e^{n}\parallel^2)+2\tau\parallel s^n\parallel^2,
\end{align}
where
$C'_2=\frac{3C^2_1+|\gamma_2|+2}{2}$.

Next, computing the discrete inner product of (\ref{truncation}) with $\Delta^{\alpha}_h e^{\bar{n}}$ and taking the real part of the resulting equation, we obtain
\begin{align}\label{estimate60}
\frac{\parallel\Lambda^{\alpha} e^{n+1}\parallel^2-\parallel\Lambda^{\alpha} e^{n-1}\parallel^2}{4\tau}+\upsilon_1\parallel \Delta^{\alpha}_h e^{\bar{n}} \parallel^2
=Re\left[-\left(P^n,\Delta^{\alpha}_h e^{\bar{n}}\right)\right]+\gamma_1\parallel\Lambda^{\alpha} e^{\bar{n}}\parallel^2+Re\left(r^n,\Delta^{\alpha}_h e^{\bar{n}}\right),
\end{align}
where Lemma~\ref{Lemma4} is used.

From (\ref{estimate49}), one  obtains
\begin{align}\label{estimate61}
\parallel P^n\parallel^2&\leq C^2_1 (\parallel e^{n}\parallel+\parallel e^{\bar{n}}\parallel+\parallel \xi^{n}\parallel)^2\nonumber\\
&\leq 3C^2_1 (\parallel e^{n-1}\parallel^2+\parallel e^{n}\parallel^2+\parallel e^{n+1}\parallel^2+\parallel \xi^{n}\parallel^2),
\end{align}
for $\tau<\tau_1, h<h_1$.

Thus, we have
\begin{align}\label{estimate62}
\left|Re\left[-\left(P^n,\Delta^{\alpha}_h e^{\bar{n}}\right)\right]\right|&\leq \parallel P^n\parallel\ \parallel \Delta^{\alpha}_h e^{\bar{n}}\parallel\nonumber\\
&\leq  (\frac{1}{2\upsilon_1}\parallel P^n\parallel^2+\frac{\upsilon_1}{2}\parallel \Delta^{\alpha}_h e^{\bar{n}}\parallel^2)\nonumber\\
&\leq \frac{3C^2_1}{2\upsilon_1}(\parallel e^{n-1}\parallel^2+\parallel e^{n}\parallel^2+\parallel e^{n+1}\parallel^2+\parallel \xi^{n}\parallel^2)+\frac{\upsilon_1}{2}\parallel \Delta^{\alpha}_h e^{\bar{n}}\parallel^2,
\end{align}
and
\begin{align}\label{estimate63}
\left| Re\left(r^n,\Delta^{\alpha}_h e^{\bar{n}}\right)\right|\leq \parallel r^n\parallel\ \parallel\Delta^{\alpha}_h e^{\bar{n}}\parallel\leq \frac{1}{2\upsilon_1}\parallel r^n\parallel^2+\frac{\upsilon_1}{2}\parallel\Delta^{\alpha}_h e^{\bar{n}}\parallel^2.
\end{align}

Substituting (\ref{estimate60}) and (\ref{estimate62}) into (\ref{estimate63}) yields
\begin{align}\label{estimate64}
&\frac{\parallel\Lambda^{\alpha} e^{n+1}\parallel^2-\parallel\Lambda^{\alpha} e^{n-1}\parallel^2}{4\tau}+\upsilon_1\parallel \Delta^{\alpha}_h e^{\bar{n}} \parallel^2\nonumber\\
\leq& \frac{3C^2_1}{2\upsilon_1}(\parallel e^{n-1}\parallel^2+\parallel e^{n}\parallel^2+\parallel e^{n+1}\parallel^2+\parallel \xi^{n}\parallel^2)+\frac{\upsilon_1}{2}\parallel \Delta^{\alpha}_h e^{\bar{n}}\parallel^2+\gamma_1\parallel\Lambda^{\alpha} e^{\bar{n}}\parallel^2+\frac{1}{2\upsilon_1}\parallel r^n\parallel^2+\frac{\upsilon_1}{2}\parallel\Delta^{\alpha}_h e^{\bar{n}}\parallel^2.
\end{align}
Thus,
\begin{align}\label{estimate65}
&\parallel\Lambda^{\alpha} e^{n+1}\parallel^2-\parallel\Lambda^{\alpha} e^{n-1}\parallel^2\nonumber\\
\leq& \frac{6C^2_1\tau}{\upsilon_1}(\parallel e^{n-1}\parallel^2+\parallel e^{n}\parallel^2+\parallel e^{n+1}\parallel^2+\parallel \xi^{n}\parallel^2)+2|\gamma_1|\tau\left(\parallel\Lambda^{\alpha} e^{n-1}\parallel^2+\parallel\Lambda^{\alpha} e^{n+1}\parallel^2\right)+\frac{2\tau}{\upsilon_1}\parallel r^n\parallel^2,
\end{align}
when $\tau<\tau_1, h< h_1$.

Similarly,  computing the discrete inner product of (\ref{truncation0}) with $\Delta_h\xi^{\bar{n}}$ and analyzing the resulting equation, one  can obtain
\begin{align}\label{estimate650}
&\parallel\Lambda^{\alpha} \xi^{n+1}\parallel^2-\parallel\Lambda^{\alpha} \xi^{n-1}\parallel^2\nonumber\\
\leq& \frac{6C_1^2\tau}{\upsilon_2}(\parallel \xi^{n-1}\parallel^2+\parallel \xi^{n}\parallel^2+\parallel \xi^{n+1}\parallel^2+\parallel e^{n}\parallel^2)+2|\gamma_2|\tau\left(\parallel\Lambda^{\alpha} \xi^{n-1}\parallel^2+\parallel\Lambda^{\alpha} \xi^{n+1}\parallel^2\right)+\frac{2\tau}{\upsilon_2}\parallel s^n\parallel^2,
\end{align}
when $\tau<\tau_1, h<h_1$.

Adding the inequalities (\ref{estimate51}), (\ref{estimate510}), (\ref{estimate65}) and (\ref{estimate650}), one obtains that
\begin{align}
&\parallel\Lambda^{\alpha} e^{n+1}\parallel^2+\parallel\Lambda^{\alpha} \xi^{n+1}\parallel^2+\parallel e^{n+1}\parallel^2+\parallel \xi^{n+1}\parallel^2-\parallel\Lambda^{\alpha} e^{n-1}\parallel-\parallel\Lambda^{\alpha} \xi^{n-1}\parallel^2-\parallel e^{n-1}\parallel^2-\parallel \xi^{n-1}\parallel^2\nonumber\\
\leq&C_3\tau(\parallel e^{n-1}\parallel^2+2\parallel e^{n}\parallel^2+\parallel e^{n+1}\parallel^2+\parallel\Lambda^{\alpha} e^{n-1}\parallel^2+2\parallel\Lambda^{\alpha} e^{n}\parallel^2+\parallel\Lambda^{\alpha} e^{n+1}\parallel^2)\nonumber\\
+&C_3\tau(\parallel \xi^{n-1}\parallel^2+2\parallel \xi^{n}\parallel^2+\parallel \xi^{n+1}\parallel^2+\parallel\Lambda^{\alpha} \xi^{n-1}\parallel^2+2\parallel\Lambda^{\alpha} \xi^{n}\parallel^2+\parallel\Lambda^{\alpha} \xi^{n+1}\parallel^2)\nonumber\\
+&\left(2+\frac{2}{\upsilon_1}\right)\tau\parallel r^n\parallel^2+\left(2+\frac{2}{\upsilon_2}\right)\tau\parallel s^n\parallel^2,
\end{align}
where $C_3=\max\left\{\frac{6C^2_1}{\upsilon_1},\frac{6C_1^2}{\upsilon_2},4C_2, 4C'_2, 2|\gamma_1|,2|\gamma_2|\right\}$.

Let $E^n=\parallel \Lambda^{\alpha}e^{n-1}\parallel^2+\parallel e^{n-1}\parallel^2+\parallel \Lambda^{\alpha} e^{n}\parallel^2+\parallel e^{n}\parallel^2+\parallel \Lambda^{\alpha}\xi^{n-1}\parallel^2+\parallel \xi^{n-1}\parallel^2+\parallel \Lambda^{\alpha} \xi^{n}\parallel^2+\parallel \xi^{n}\parallel^2$, then one has
\begin{align}
&E^{n+1}-E^{n}\leq C_3\tau(E^{n+1}+E^{n})+\left(2+\frac{2}{\upsilon_1}\right)\tau\parallel r^n\parallel^2+\left(2+\frac{2}{\upsilon_2}\right)\tau\parallel s^n\parallel^2,
\end{align}
which is equivalent to
\begin{align}
(1-C_3\tau)(E^{n+1}-E^n)\leq 2C_3\tau E^n+\left(2+\frac{2}{\upsilon_1}\right)\tau\parallel r^n\parallel^2+\left(2+\frac{2}{\upsilon_2}\right)\tau\parallel s^n\parallel^2.
\end{align}
When $\tau<\tau_2=\frac{1}{2C_3}$, then $1-C_3\tau>\frac{1}{2}$, (\ref{estimate47}) gives
\begin{align}\label{estimate67}
E^{n+1}-E^{n}\leq 4C_3\tau E^n+\left(4+\frac{4}{\upsilon_1}\right)\tau\parallel r^n\parallel^2+\left(4+\frac{4}{\upsilon_2}\right)\tau\parallel s^n\parallel^2.
\end{align}

Replacing $n$ by $k$ in (\ref{estimate67}) and summing over $k$ from $1$ to $n$ yields
\begin{align}\label{estimate68}
E^{n+1}-E^{1}\leq  4C_3\tau \sum^{n}_{k=1}E^k+\left(4+\frac{4}{\upsilon_1}\right)\tau\sum^{n}_{k=1}\parallel r^k\parallel^2+\left(4+\frac{4}{\upsilon_2}\right)\tau\sum^{n}_{k=1}\parallel s^k\parallel^2.
\end{align}

From Lemma~\ref{Lemma8}, one gets
\begin{equation}\label{rs}
  \begin{aligned}
&\tau\sum^{n}_{k=1}\parallel r^k \parallel^2\leq \tau nC_R^2(\tau^2+h^2)^2\leq C_R^2T(\tau^2+h^2)^2, \\ &\tau\sum^{n}_{k=1}\parallel s^k \parallel^2\leq \tau nC_R^2(\tau^2+h^2)^2\leq C_R^2T(\tau^2+h^2)^2.
\end{aligned}
\end{equation}
Lemma~\ref{Lemma11} yields
\begin{align*}
\parallel e^{1}\parallel^2=(e^1,e^1)\leq& h\sum^{M-1}_{j=1}|e^1_j||e^1_j|\nonumber\\
\leq&(b-a)\max_{1\leq j\leq M-1}|e^1_j|^2\nonumber\\
\leq&C_e^2(b-a)(\tau^2+\tau h^2)^2\nonumber\\
\leq&C_e^2(b-a)(\tau^2+h^2)^2,
\end{align*}
and
\begin{align*}
\parallel\Lambda^{\alpha} e^{1}\parallel^2=(\Delta^{\alpha}_he^1,e^1)\leq& h\sum^{M-1}_{j=1}|\Delta^{\alpha}_he^1_j||e^1_j|\nonumber\\
\leq&(b-a)\max_{1\leq j\leq M-1}|\Delta^{\alpha}_he^1_j|\max_{1\leq j\leq M-1}|e^1_j|\nonumber\\
\leq&C_e^2(b-a)(\tau^2+\tau h^{2+\alpha})(\tau^2+\tau h^2)\nonumber\\
\leq&C_e^2(b-a)(\tau^2+h^2)^2,
\end{align*}
for $\tau \leq 1$ and $h\leq 1.$ Thus,
\begin{align}
\parallel e^{1}\parallel\leq C_e\sqrt{b-a}(\tau^2+ h^2), \quad \parallel\Lambda^{\alpha} e^{1}\parallel\leq C_e\sqrt{b-a}(\tau^2+ h^2).
\end{align}

Similarly,
\begin{align}
\parallel \xi^{1}\parallel\leq C_e\sqrt{b-a}(\tau^2+ h^2), \quad \parallel\Lambda^{\alpha} \xi^{1}\parallel\leq C_e\sqrt{b-a}(\tau^2+ h^2).
\end{align}
In addition, $\parallel e^0\parallel=0,\parallel \Lambda^{\alpha} e^0\parallel=0,$ $\parallel \xi^0\parallel=0,\parallel \Lambda^{\alpha} \xi^0\parallel=0$.
Thus, \begin{align}
  E^1&=\parallel\Lambda^{\alpha} e^{1}\parallel^2+\parallel e^{1}\parallel^2+\parallel\Lambda^{\alpha} e^{0}\parallel^2+\parallel e^{0}\parallel^2+\parallel\Lambda^{\alpha} \xi^{1}\parallel^2+\parallel \xi^{1}\parallel^2+\parallel\Lambda^{\alpha} \xi^{0}\parallel^2+\parallel \xi^{0}\parallel^2 \notag \\
  & \leq 4C_e^2(b-a)(\tau^2+ h^2)^2.\label{e1est}
\end{align}

Substituting (\ref{rs}) and (\ref{e1est}) into (\ref{estimate68}) gives
\begin{align}\label{estimate70}
E^{n+1}\leq 4C_3\tau \sum^{n}_{k=1}E^k+C_4(\tau^2+h^2)^2,
\end{align}
where $C_4 = 4C_e^2(b-a)+(8+\frac{4}{\upsilon_1}+\frac{4}{\upsilon_2})C_R^2T$.

By using Lemma~\ref{Lemma6}, one has
\begin{align}
E^{n+1}\leq C_4(\tau^2+h^2)^2 e^{4C_3\tau n}\leq C_4 e^{4C_3T}(\tau^2+h^2)^2.
\end{align}
From Lemma~\ref{Lemma3} and Lemma~\ref{Lemma4}, one gets
\begin{align}
&\parallel e^{n+1}\parallel^2+\left(\frac{2}{\pi}\right)^{\alpha}|e^{n+1}|^2_{H^{\frac{\alpha}{2}}}\leq \parallel e^{n+1}\parallel^2+\parallel \Lambda^{\alpha} e^{n+1}\parallel^2 \leq  C_4 e^{4C_3T}(\tau^2+h^2)^2, \\
&\parallel \xi^{n+1}\parallel^2+\left(\frac{2}{\pi}\right)^{\alpha}|\xi^{n+1}|^2_{H^{\frac{\alpha}{2}}}\leq \parallel \xi^{n+1}\parallel^2+\parallel \Lambda^{\alpha} \xi^{n+1}\parallel^2 \leq  C_4 e^{4C_3T}(\tau^2+h^2)^2.
\end{align}
Thus, one obtains that
\begin{align}
&\parallel e^{n+1}\parallel^2_{H^{\frac{\alpha}{2}}}\leq C_4 e^{4C_3T}\left(\frac{\pi}{2}\right)^{\alpha}(\tau^2+h^2)^2, \\
&\parallel \xi^{n+1}\parallel^2_{H^{\frac{\alpha}{2}}}\leq C_4 e^{4C_3T}\left(\frac{\pi}{2}\right)^{\alpha}(\tau^2+h^2)^2.
\end{align}

Therefore,
\begin{align}
\parallel e^{n+1}\parallel_{H^{\frac{\alpha}{2}}}\leq C_5(\tau^2+h^2),\quad \parallel \xi^{n+1}\parallel_{H^{\frac{\alpha}{2}}}\leq C_5(\tau^2+h^2)
\end{align}
where $C_5=\sqrt{C_4 e^{4C_3T}\left(\frac{\pi}{2}\right)^{\alpha}}$.

By Lemma~\ref{Lemma2}, we obtain
\begin{align}\label{estimate71}
\parallel e^{n+1}\parallel_{\infty}\leq C_\alpha \parallel e^{n+1}\parallel_{H^{\frac{\alpha}{2}}} \leq C_\alpha C_5(\tau^2+h^2),\quad \parallel \xi^{n+1}\parallel_{\infty}\leq C_\alpha \parallel \xi^{n+1}\parallel_{H^{\frac{\alpha}{2}}} \leq C_\alpha C_5(\tau^2+h^2).
\end{align}

Now we take $C_0=\max\left\{C_\alpha C_5, C_e\right\}$. Once $C_0$ is fixed, the condition for $\tau_1, h_1$, i.e., $\tau_1^2+h_1^2<1/C_0$ can be used to determine $\tau_1, h_1$.

Thus, let $\tau_0=\min\left\{\tau_1,\tau_2, 1\right\}$ and $h_0=\min\left\{h_1, 1\right\}$, then (\ref{estimate1}) is valid for $n+1$.
The induction is closed. This completes the proof. \qed

\begin{theorem}\label{theorem5}
Suppose that the solution of problem (\ref{IVP1})-(\ref{IVP3}) is smooth enough, the solution $U^n$ of the difference scheme (\ref{diff1})-(\ref{diff3})  is bounded in the $L^{\infty}$-norm for $\tau<\tau_0$ and $h<h_0$, i.e.,
 \begin{align}
\parallel U^n\parallel_{\infty}\leq C_*,\quad \parallel V^n\parallel_{\infty}\leq C_* \quad  1\leq n \leq N.
\end{align}
where $\tau_0, h_0$ are the same positive constants in Theorem~\ref{theorem3}.
 \end{theorem}

 \noindent {\bf Proof.} From Theorem~\ref{theorem3}, one has
 \begin{align}
\parallel U^n\parallel_{\infty}&\leq \parallel u^n\parallel_{\infty}+\parallel e^n\parallel_{\infty}, \nonumber\\
&\leq C_m+C_0(\tau^2+h^2), \quad  1\leq n \leq N,
\end{align}
and
 \begin{align}
\parallel V^n\parallel_{\infty}&\leq \parallel v^n\parallel_{\infty}+\parallel \xi^n\parallel_{\infty}, \nonumber\\
&\leq C_m+C_0(\tau^2+h^2), \quad  1\leq n \leq N,
\end{align}
for $\tau<\tau_0,h<h_0$.

Since $C_0(\tau^2+h^2)<1$ for  $\tau<\tau_0, h<h_0$, one has
 \begin{align}
\parallel U^n\parallel_{\infty}\leq  C_*,\quad \parallel V^n\parallel_{\infty}\leq  C_*, \quad  2\leq n \leq N,
\end{align}
where $C_*=C_m+1$. This completes the proof.  \qed
%
%
%

\subsection{Existence and uniqueness}
 \begin{theorem}\label{theorem2}
Suppose that the solution of problem (\ref{IVP1})-(\ref{IVP3}) is smooth enough. Then the difference scheme (\ref{diff1})-(\ref{diff3}) is uniquely solvable for $\tau <\min\left\{\tau_0,\frac{1}{|\gamma_1|+(|\kappa_1|+|\delta_1|)|C_*|^2},\frac{1}{|\gamma_2|+(|\kappa_1|+|\delta_1|)|C_*|^2}\right\}$ and $h<h_0$, where $\tau_0, h_0$ are the same positive constants in Theorem~\ref{theorem3}.
 \end{theorem}
 {\bf Proof.} To prove the theorem, we proceed by the mathematical induction. Obviously,  $U^1$ and $V^1$ can be uniquely determined by (\ref{diffinitialaa}) and (\ref{diffinitialaabb}).
 Suppose  $U^1,U^2,\cdots,U^{n}, V^1,V^2,\cdots, V^{n} (1\leq n\leq N-1)$ are obtained uniquely, we now show that $U^{n+1}$ and $V^{n+1}$ are uniquely determined by (\ref{diff1})-(\ref{diff1l}).

 Assume that  $U^{n+1,1},U^{n+1,2}$ are two solutions of (\ref{diff1}) and let $W^{n+1}=U^{n+1,1}-U^{n+1,2}$,  then it is easy to verify that  $W^{n+1}$ satisfies the following equation:
 \begin{align}\label{unique1}
\frac{1}{2\tau}W^{n+1}_j+\frac{1}{2}(\upsilon_1+i\eta_1)\Delta^{\alpha}_hW^{n+1}_{j}+\frac{1}{2}\left((\kappa_1+i\zeta_1)|U^{n}_{j}|^2+(\delta_1+i\beta_1)|V^{n}_{j}|^2\right)W^{n+1}_{j}-\frac{1}{2}\gamma_1 W^{n+1}_{j}=0,\quad 0<j<M.
 \end{align}

Computing the inner product of (\ref{unique1}) with $W^{n+1}$ and taking the real part of the resulting equation, we have
  \begin{align}
 (\frac{1}{2\tau}-\frac{\gamma_1}{2})\parallel W^{n+1}\parallel^2+\frac{\upsilon_1}{2}\parallel \Lambda^{\alpha}W^{n+1}\parallel^2+\frac{\kappa_1}{2}{\parallel U^nW^{n+1}\parallel^2}+\frac{\delta_1}{2}{\parallel V^nW^{n+1}\parallel^2}=0,
 \end{align}
 where Lemma~\ref{Lemma4} is used.

If $\tau<\tau_0, h<h_0$, then
 \begin{align}
 0=&(\frac{1}{2\tau}-\frac{\gamma_1}{2})\parallel W^{n+1}\parallel^2+\frac{\upsilon_1}{2}\parallel \Lambda^{\alpha}W^{n+1}\parallel^2+\frac{\kappa_1}{2}{\parallel U^nW^{n+1}\parallel^2}+\frac{\delta_1}{2}{\parallel V^nW^{n+1}\parallel^2}\nonumber\\
 \geq &(\frac{1}{2\tau}-\frac{|\gamma_1|}{2}-\frac{|\kappa_1|}{2}C_*^2-\frac{|\delta_1|}{2}C_*^2)\parallel W^{n+1}\parallel^2+\frac{\upsilon_1}{2}\parallel \Lambda^{\alpha}W^{n+1}\parallel^2, \label{unique2}
 \end{align}
where Theorem~\ref{theorem5} is used.

 If $\tau < \min\left\{\frac{1}{|\gamma_1|+(|\kappa_1|+|\delta_1|)|C_*|^2},\frac{1}{|\gamma_2|+(|\kappa_1|+|\delta_1|)|C_*|^2}\right\}$, then $\frac{1}{2\tau}-\frac{\gamma_1}{2}-\frac{|\kappa_1|}{2}C_*^2-\frac{|\delta_1|}{2}C_*^2> 0$.
  Since $\upsilon_1>0$, (\ref{unique2}) implies
  \begin{align}\label{unique3}
\parallel W^{n+1}\parallel=\parallel \Lambda^{\alpha}W^{n+1}\parallel=0.
 \end{align}

 Therefore, (\ref{unique1})  has only a trivial solution.  This proves the uniqueness of the numerical solution $U^{n+1}$. The proof of the uniqueness for $V^{n+1}$ is similar. This completes the proof of the theorem. \qed

\begin{remark}
  The method presented in the above is second-order accurate both in time and space variables. And the second-order scheme (\ref{diff1}) is easy to extend to the case of spatial fourth-order. Let $\mathcal{A}^{\alpha}_{x}$ {be the average operator defined as~\cite{Hao2016}
  \begin{equation}
    \mathcal{A}^{\alpha}_{x} U_j= \frac{\alpha}{24} U_{j-1} + \left(1-\frac{\alpha}{12}\right) U_{j} + \frac{\alpha}{24} U_{j+1},
      \end{equation}}
  then our spatial fourth-order method is stated as follows,
\begin{align}
\mathcal{A}^{\alpha}_{x}\delta_{t} U^{n}_j+(\upsilon_1+i\eta_1)\Delta^{\alpha}_hU^{\bar{n}}_{j}+\mathcal{A}^{\alpha}_{x}\left((\kappa_1+i\zeta_1)|U^{n}_{j}|^2+(\delta_1+i\beta_1)|V^{n}_{j}|^2\right)U^{\bar{n}}_{j}-\gamma_1 \mathcal{A}^{\alpha}_{x} U^{\bar{n}}_{j}&=0,\quad 0<j<M,\ 1<n<N,\label{fourth1}\\
\mathcal{A}^{\alpha}_{x}\delta_{t} V^{n}_j+(\upsilon_2+i\eta_2)\Delta^{\alpha}_hV^{\bar{n}}_{j}+\mathcal{A}^{\alpha}_{x}\left((\kappa_2+i\zeta_2)|U^{n}_{j}|^2+(\delta_2+i\beta_2)|V^{n}_{j}|^2\right)V^{\bar{n}}_{j}-\gamma_2 \mathcal{A}^{\alpha}_{x} V^{\bar{n}}_{j}&=0,\quad 0<j<M,\ 1<n<N,\label{fourth2}
\end{align}
where the method is second-order accurate in time and fourth-order accurate in space.
Replacing the average operator $\mathcal{A}^{\alpha}_{x}$ by the identity operator in (\ref{fourth1}) and (\ref{fourth2}), yields the second-order scheme (\ref{diff1}). To start the method, $U^1$ can also be computed through (\ref{diffinitialaa}), where the truncation error  is
\begin{equation*}
U^1_j-u^1_j=O(\tau h^2)+O(\tau^2)\leq O(\tau^2)+O(h^4).
\end{equation*}

Following the proof in the previous sections, one can also prove that  the above method (\ref{fourth1})-(\ref{fourth2}) is also unconditionally stable, second-order accurate in time and fourth-order accurate in space, and a pointwise error estimate can be achieved.

In order to increase the time accuracy, the following Richardson extrapolation for the final step numerical solution is used:
\begin{equation}\label{Richardson}
\begin{aligned}
    \widetilde{U}^N(\Delta t, h) &= \frac{4}{3}U^N(\Delta t, h) - \frac{1}{3} U^{N/2}(2\Delta t, h), \\
    \widetilde{V}^N(\Delta t, h) &= \frac{4}{3}V^N(\Delta t, h) - \frac{1}{3} V^{N/2}(2\Delta t, h),
\end{aligned}
\end{equation}
where $U^N(\Delta t, h)$, $U^{{N}/{2}}(2\Delta t, h)$ are numerical solutions at the final step by using spatial meshsizes $h$ and time step $\Delta t$, $2\Delta t$, respectively.

\end{remark}

\section{Numerical results}\label{section5}
    In this section, we present some numerical results of the proposed difference scheme (\ref{diff1})-(\ref{diff3}) to support our theoretical findings.

\begin{example}\label{ex1}
In order to test the accuracy of the proposed scheme, we consider the following system with source terms:
\begin{align*}
&u_t+(1+i)(-\Delta)^{\frac{\alpha}{2}}u+\left((-1-i)|u|^2+(1+i)|v|^2\right)u- u=f(x ,t),\quad x\in (0,1), \quad 0<t<1,\\
&v_t+(1-i)(-\Delta)^{\frac{\alpha}{2}}v+\left((1+i)|u|^2+(1-i)|v|^2\right)v+v=g(x, t), \quad x\in (0,1), \quad 0<t<1,
\end{align*}
with homogeneous boundary conditions
\begin{equation*}
  u(0, t) = u(1, t) =0,\quad v(0, t) =v(1,t) =0.
\end{equation*}
The initial conditions, and the source terms $f(x, t)$ and $g(x, t)$ are determined by the exact solutions
\begin{equation*}
  u(x, t) = \exp(-t) x^4(1-x)^4,\quad v(x, t) = (t+1)^3 x^4(1-x)^4.
\end{equation*}
\end{example}

Table~\ref{tab1} and Table~\ref{tab2} list the errors and the convergence orders for the method (\ref{diff1})-(\ref{diff3}) with $\alpha = 1.2, 1.5, 1.8, 2.0$ in the $L^\infty$-norm, respectively. As we can see that these results confirm the  second-order convergence both in time and space  variables.  Table~\ref{tab3} and Table~\ref{tab4} list the errors and the convergence orders for the method (\ref{fourth1})-(\ref{Richardson}) with $\alpha = 1.2, 1.5, 1.8, 2.0$ in the $L^\infty$-norm, respectively. As we can see that these results show that the method (\ref{fourth1})-(\ref{Richardson}) is  fourth-order convergence both in time and space variables.

\begin{figure}[!tbp]
\centering
\subfigure[Numerical error of $U$ at $T=1$]{
\includegraphics[width=0.45\textwidth]{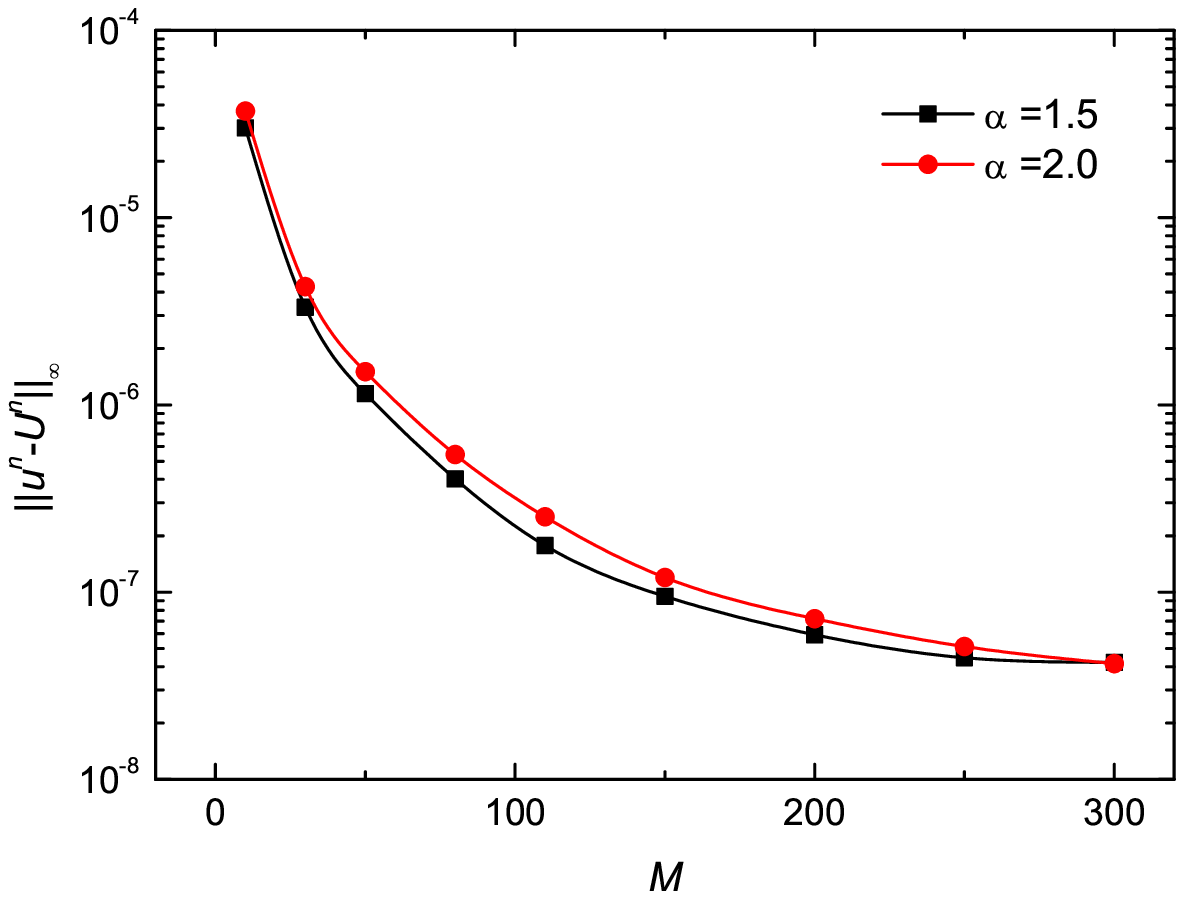}}
\subfigure[Numerical error of $V$ at $T=1$]{
\includegraphics[width=0.45\textwidth]{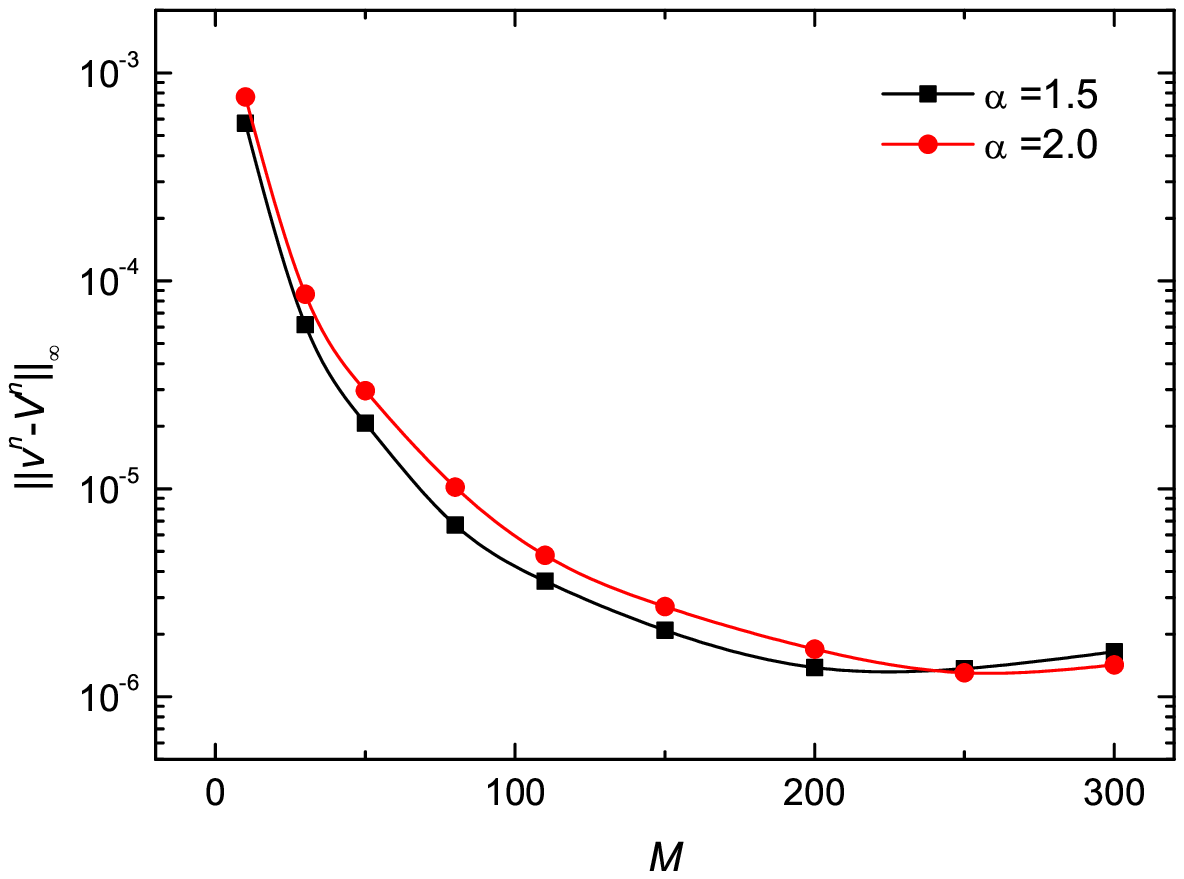}} \\
\caption{Numerical errors obtained by the second-order method for fixed $\tau=0.01$.}
\label{Fig1}
\end{figure}

\begin{figure}[!tbp]
\centering
\subfigure[Numerical error of $\widetilde{U}$ at $T=1$]{
\includegraphics[width=0.45\textwidth]{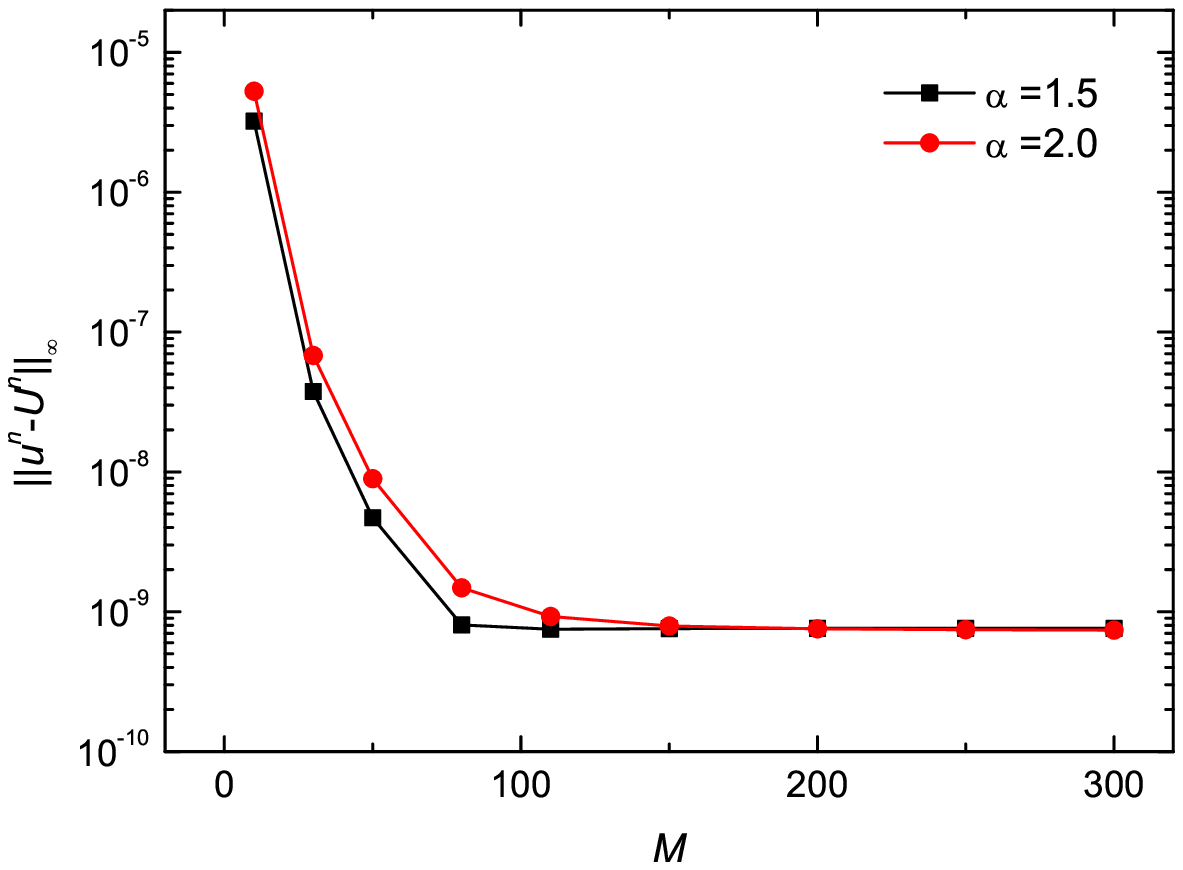}}
\subfigure[Numerical error of $\widetilde{V}$ at $T=1$]{
\includegraphics[width=0.45\textwidth]{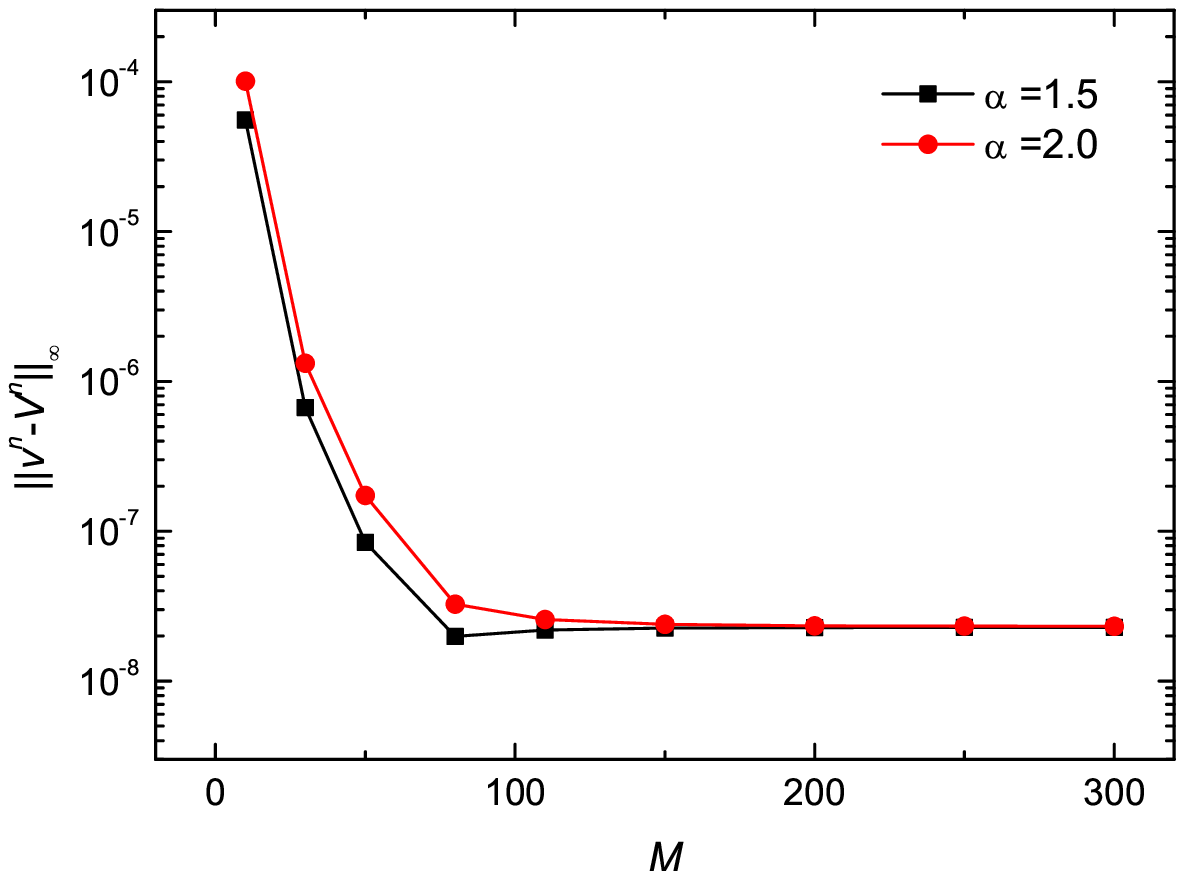}} \\
\caption{Numerical errors obtained by the fourth-order method for fixed $\tau=0.001$.}
\label{Fig2}
\end{figure}

In order the illustrate the unconditionally stability of the our methods, we fix $\tau$ and vary $h$, results for $\alpha=1.5$ and $\alpha=2$ are plotted in Figure~\ref{Fig1} and Figure~\ref{Fig2}. As one can see that these results clearly show that the time step is not related to the spatial meshsize, and as the spatial meshsize goes to zero, the dominant error comes from the temporal part.

\begin{table}[!tbp]
\tabcolsep=5pt
\caption{$L^\infty$-norm errors and their convergence orders of ${U}$ obtained by the second-order scheme for Example \ref{ex1}.}\label{tab1}
 \centering
  \begin{tabular}{lllllllllllll}
    \hline
   \multirow{2}*{$\tau$} & \multirow{2}*{$h$}  & \multicolumn{2}{l}{$\alpha=1.2$} &  & \multicolumn{2}{l}{$\alpha=1.5$} & & \multicolumn{2}{l}{$\alpha=1.8$} & & \multicolumn{2}{l}{$\alpha=2.0$}  \\
    \cline{3-4} \cline{6-7} \cline{9-10} \cline{12-13}&  &   $\parallel u^n-{U}^n\parallel_\infty$ &  order  & &  $\parallel u^n-{U}^n\parallel_\infty$ &  order&      &$\parallel u^n-{U}^n\parallel_\infty$& order  & &    $\parallel u^n-{U}^n\parallel_\infty$  &order   \\
\hline
$1/32$  &$1/32$ &   1.72e$-$06   &   -    &&   2.25e$-$06  &    -      &&  2.76e$-$06  &    -    &&   3.14e$-$06   &   -           \\
$1/64$  &$1/64$ &   4.31e$-$07   &   2.00 &&   5.63e$-$07  &    2.00   &&  6.95e$-$07  &    1.99 &&   7.85e$-$07   &   2.00        \\
$1/128$ &$1/128$&   1.08e$-$07   &   2.00 &&   1.41e$-$07  &    2.00   &&  1.74e$-$07  &    2.00 &&   1.96e$-$07   &   2.00        \\
$1/256$ &$1/256$&   2.69e$-$08   &   2.00 &&   3.52e$-$08  &    2.00   &&  4.35e$-$08  &    2.00 &&   4.91e$-$08   &   2.00        \\
$1/512$ &$1/512$&   6.72e$-$09   &   2.00 &&   8.80e$-$09  &    2.00   &&  1.09e$-$08  &    2.00 &&   1.23e$-$08   &   2.00        \\
\hline
\end{tabular}
\end{table}

\begin{table}[!tbp]
\tabcolsep=5pt
\caption{$L^\infty$-norm errors and their convergence orders of ${V}$ obtained by the second-order scheme for Example \ref{ex1}.}\label{tab2}
 \centering
  \begin{tabular}{lllllllllllll}
    \hline
   \multirow{2}*{$\tau$} & \multirow{2}*{$h$}  & \multicolumn{2}{l}{$\alpha=1.2$} &  & \multicolumn{2}{l}{$\alpha=1.5$} & & \multicolumn{2}{l}{$\alpha=1.8$} & & \multicolumn{2}{l}{$\alpha=2.0$}  \\
    \cline{3-4} \cline{6-7} \cline{9-10} \cline{12-13}&  &   $\parallel v^n-{V}^n\parallel_\infty$ &  order  & &  $\parallel v^n-{V}^n\parallel_\infty$ &  order&      &$\parallel v^n-{V}^n\parallel_\infty$& order  & &    $\parallel v^n-{V}^n\parallel_\infty$  &order   \\
\hline
$1/32$  &$1/32$     &  3.32e$-$05    &  -       &&    4.23e$-$05  &    -    &&    5.11e$-$05   &   -      &&    5.67e$-$05   &   -              \\
$1/64$  &$1/64$     &  8.27e$-$06    &  2.01    &&    1.05e$-$05  &    2.01 &&    1.26e$-$05   &   2.01   &&    1.40e$-$05   &   2.02           \\
$1/128$ &$1/128$    &  2.06e$-$06    &  2.00    &&    2.62e$-$06  &    2.00 &&    3.16e$-$06   &   2.00   &&    3.50e$-$06   &   2.00           \\
$1/256$ &$1/256$    &  5.16e$-$07    &  2.00    &&    6.56e$-$07  &    2.00 &&    7.89e$-$07   &   2.00   &&    8.75e$-$07   &   2.00           \\
$1/512$ &$1/512$    &  1.29e$-$07    &  2.00    &&    1.64e$-$07  &    2.00 &&    1.97e$-$07   &   2.00   &&    2.19e$-$07   &   2.00           \\
\hline
\end{tabular}
\end{table}

\begin{table}[!tbp]
\tabcolsep=5pt
\caption{$L^\infty$-norm errors and their convergence orders of $\widetilde{U}$ obtained by the fourth-order scheme for Example \ref{ex1}.}\label{tab3}
 \centering
  \begin{tabular}{lllllllllllll}
    \hline
   \multirow{2}*{$\tau$} & \multirow{2}*{$h$}  & \multicolumn{2}{l}{$\alpha=1.2$} &  & \multicolumn{2}{l}{$\alpha=1.5$} & & \multicolumn{2}{l}{$\alpha=1.8$} & & \multicolumn{2}{l}{$\alpha=2.0$}  \\
    \cline{3-4} \cline{6-7} \cline{9-10} \cline{12-13}&  &   $\parallel u^n-\widetilde{U}^n\parallel_\infty$ &  order  & &  $\parallel u^n-\widetilde{U}^n\parallel_\infty$ &  order&      &$\parallel u^n-\widetilde{U}^n\parallel_\infty$& order  & &    $\parallel u^n-\widetilde{U}^n\parallel_\infty$  &order   \\
\hline
$1/32$  &$1/32$ &   1.61e$-$08  &    -     &&  2.83e$-$08  &    -    &&    1.05e$-$07  &    -      && 1.68e$-$07  &    -      \\
$1/64$  &$1/64$ &   8.78e$-$10  &    4.20  &&  1.81e$-$09  &    3.97 &&    2.66e$-$09  &    5.31   && 3.15e$-$09  &    5.74      \\
$1/128$ &$1/128$&   5.75e$-$11  &    3.93  &&  1.05e$-$10  &    4.11 &&    1.48e$-$10  &    4.16   && 2.24e$-$10  &    3.81      \\
$1/256$ &$1/256$&   3.93e$-$12  &    3.87  &&  6.14e$-$12  &    4.09 &&    8.99e$-$12  &    4.05   && 1.30e$-$11  &    4.11      \\
$1/512$ &$1/512$&   2.67e$-$13  &    3.88  &&  3.63e$-$13  &    4.08 &&    5.40e$-$13  &    4.06   && 8.07e$-$13  &    4.01      \\
\hline
\end{tabular}
\end{table}

\begin{table}[!tbp]
\tabcolsep=5pt
\caption{$L^\infty$-norm errors and their convergence orders of $\widetilde{V}$ obtained by the fourth-order scheme for Example \ref{ex1}.}\label{tab4}
 \centering
  \begin{tabular}{lllllllllllll}
    \hline
   \multirow{2}*{$\tau$} & \multirow{2}*{$h$}  & \multicolumn{2}{l}{$\alpha=1.2$} &  & \multicolumn{2}{l}{$\alpha=1.5$} & & \multicolumn{2}{l}{$\alpha=1.8$} & & \multicolumn{2}{l}{$\alpha=2.0$}  \\
    \cline{3-4} \cline{6-7} \cline{9-10} \cline{12-13}&  &   $\parallel v^n-\widetilde{V}^n\parallel_\infty$ &  order  & &  $\parallel v^n-\widetilde{V}^n\parallel_\infty$ &  order&      &$\parallel v^n-\widetilde{V}^n\parallel_\infty$& order  & &    $\parallel v^n-\widetilde{V}^n\parallel_\infty$  &order   \\
\hline
$1/32$  &$1/32$     &  3.42e$-$07  &    -      &&  5.19e$-$07  &    -      &&   1.08e$-$06 &     -      &&  1.29e$-$06   &   -          \\
$1/64$  &$1/64$     &  1.96e$-$08  &    4.12   &&  3.04e$-$08  &    4.09   &&   4.62e$-$08 &     4.55   &&  5.69e$-$08   &   4.50       \\
$1/128$ &$1/128$    &  1.16e$-$09  &    4.08   &&  1.78e$-$09  &    4.09   &&   2.79e$-$09 &     4.05   &&  4.10e$-$09   &   3.80       \\
$1/256$ &$1/256$    &  7.76e$-$11  &    3.90   &&  1.05e$-$10  &    4.08   &&   1.68e$-$10 &     4.05   &&  2.49e$-$10   &   4.04       \\
$1/512$ &$1/512$    &  5.62e$-$12  &    3.79   &&  6.23e$-$12  &    4.08   &&   1.01e$-$11 &     4.05   &&  1.56e$-$11   &   4.00       \\
\hline
\end{tabular}
\end{table}

\begin{table}[!tbp]
\tabcolsep=5pt
\caption{$L^\infty$-norm errors and their convergence orders of ${U}$ obtained by the second-order scheme for Example \ref{ex2}.}\label{tab5}
 \centering
  \begin{tabular}{lllllllllllll}
    \hline
   \multirow{2}*{$\tau$} & \multirow{2}*{$h$}  & \multicolumn{2}{l}{$\alpha=1.2$} &  & \multicolumn{2}{l}{$\alpha=1.5$} & & \multicolumn{2}{l}{$\alpha=1.8$} & & \multicolumn{2}{l}{$\alpha=2.0$}  \\
    \cline{3-4} \cline{6-7} \cline{9-10} \cline{12-13}&  &   $\parallel u^n-{U}^n\parallel_\infty$ &  order  & &  $\parallel u^n-{U}^n\parallel_\infty$ &  order&      &$\parallel u^n-{U}^n\parallel_\infty$& order  & &    $\parallel u^n-{U}^n\parallel_\infty$  &order   \\
\hline
$1/8$  &$1/8$     & 1.05e$-$01  &    -    &&   2.24e$-$02    &  -    &&  4.13e$-$02   &   -      &&   4.66e$-$02   &   -            \\
$1/16$ &$1/16$    & 3.01e$-$02  &    1.80 &&   6.59e$-$03    &  1.77 &&  1.21e$-$02   &   1.77   &&   1.27e$-$02   &   1.87         \\
$1/32$ &$1/32$    & 7.80e$-$03  &    1.95 &&   1.76e$-$03    &  1.90 &&  3.12e$-$03   &   1.96   &&   3.24e$-$03   &   1.97         \\
$1/64$ &$1/64$    & 1.99e$-$03  &    1.97 &&   4.43e$-$04    &  1.99 &&  7.69e$-$04   &   2.02   &&   8.03e$-$04   &   2.01         \\
\hline
\end{tabular}
\end{table}

\begin{table}[!tbp]
\tabcolsep=5pt
\caption{$L^\infty$-norm errors and their convergence orders of ${V}$ obtained by the second-order scheme for Example \ref{ex2}.}\label{tab6}
 \centering
  \begin{tabular}{lllllllllllll}
    \hline
   \multirow{2}*{$\tau$} & \multirow{2}*{$h$}  & \multicolumn{2}{l}{$\alpha=1.2$} &  & \multicolumn{2}{l}{$\alpha=1.5$} & & \multicolumn{2}{l}{$\alpha=1.8$} & & \multicolumn{2}{l}{$\alpha=2.0$}  \\
    \cline{3-4} \cline{6-7} \cline{9-10} \cline{12-13}&  &   $\parallel v^n-{V}^n\parallel_\infty$ &  order  & &  $\parallel v^n-{V}^n\parallel_\infty$ &  order&      &$\parallel v^n-{V}^n\parallel_\infty$& order  & &    $\parallel v^n-{V}^n\parallel_\infty$  &order   \\
\hline
$1/8$  &$1/8$     & 8.15e$-$02   &   -     &&   2.80e$-$02   &   -       &&   4.91e$-$02  &    -     &&  5.20e$-$02   &   -             \\
$1/16$ &$1/16$    & 2.37e$-$02   &   1.78  &&   7.59e$-$03   &   1.88    &&   1.40e$-$02  &    1.81  &&  1.41e$-$02   &   1.89          \\
$1/32$ &$1/32$    & 6.17e$-$03   &   1.94  &&   1.94e$-$03   &   1.97    &&   3.59e$-$03  &    1.96  &&  3.57e$-$03   &   1.98          \\
$1/64$ &$1/64$    & 1.58e$-$03   &   1.97  &&   4.85e$-$04   &   2.00    &&   8.88e$-$04  &    2.02  &&  8.85e$-$04   &   2.01          \\
\hline
\end{tabular}
\end{table}

 \begin{table}[!tbp]
\tabcolsep=5pt
\caption{$L^\infty$-norm errors and their convergence orders of $\widetilde{U}$ obtained by the fourth-order scheme for Example \ref{ex2}.}\label{tab7}
 \centering
  \begin{tabular}{lllllllllllll}
    \hline
   \multirow{2}*{$\tau$} & \multirow{2}*{$h$}  & \multicolumn{2}{l}{$\alpha=1.2$} &  & \multicolumn{2}{l}{$\alpha=1.5$} & & \multicolumn{2}{l}{$\alpha=1.8$} & & \multicolumn{2}{l}{$\alpha=2.0$}  \\
    \cline{3-4} \cline{6-7} \cline{9-10} \cline{12-13}&  &   $\parallel u^n-\widetilde{U}^n\parallel_\infty$ &  order  & &  $\parallel u^n-\widetilde{U}^n\parallel_\infty$ &  order&      &$\parallel u^n-\widetilde{U}^n\parallel_\infty$& order  & &    $\parallel u^n-\widetilde{U}^n\parallel_\infty$  &order   \\
\hline
$1/8$  &$1/8$     &  5.65e$-$02   &   -     && 2.22e$-$02  &    -     &&   2.58e$-$02    &  -     &&   2.15e$-$02  &    -              \\
$1/16$ &$1/16$    &  8.41e$-$03   &   2.75  && 3.32e$-$03  &    2.74  &&   2.17e$-$03    &  3.57  &&   1.46e$-$03  &    3.88           \\
$1/32$ &$1/32$    &  7.56e$-$04   &   3.47  && 2.86e$-$04  &    3.54  &&   1.43e$-$04    &  3.92  &&   9.45e$-$05  &    3.95           \\
$1/64$ &$1/64$    &  5.35e$-$05   &   3.82  && 1.93e$-$05  &    3.89  &&   9.18e$-$06    &  3.96  &&   6.11e$-$06  &    3.95           \\
\hline
\end{tabular}
\end{table}

\begin{table}[!tbp]
\tabcolsep=5pt
\caption{$L^\infty$-norm errors and their convergence orders of $\widetilde{V}$ obtained by the fourth-order scheme for Example \ref{ex2}.}\label{tab8}
 \centering
  \begin{tabular}{lllllllllllll}
    \hline
   \multirow{2}*{$\tau$} & \multirow{2}*{$h$}  & \multicolumn{2}{l}{$\alpha=1.2$} &  & \multicolumn{2}{l}{$\alpha=1.5$} & & \multicolumn{2}{l}{$\alpha=1.8$} & & \multicolumn{2}{l}{$\alpha=2.0$}  \\
    \cline{3-4} \cline{6-7} \cline{9-10} \cline{12-13}&  &   $\parallel v^n-\widetilde{V}^n\parallel_\infty$ &  order  & &  $\parallel v^n-\widetilde{V}^n\parallel_\infty$ &  order&      &$\parallel v^n-\widetilde{V}^n\parallel_\infty$& order  & &    $\parallel v^n-\widetilde{V}^n\parallel_\infty$  &order   \\
\hline
$1/8$  &$1/8$     &   4.91e$-$02   &   1.55  &&   2.13e$-$02  &    2.45   &&    2.62e$-$02   &   2.41  &&   2.14e$-$02  &    2.86           \\
$1/16$ &$1/16$    &   7.69e$-$03   &   2.67  &&   3.27e$-$03  &    2.71   &&    2.14e$-$03   &   3.62  &&   1.41e$-$03  &    3.93           \\
$1/32$ &$1/32$    &   7.16e$-$04   &   3.43  &&   2.83e$-$04  &    3.53   &&    1.40e$-$04   &   3.93  &&   9.17e$-$05  &    3.94           \\
$1/64$ &$1/64$    &   5.15e$-$05   &   3.80  &&   1.92e$-$05  &    3.88   &&    8.98e$-$06   &   3.97  &&   5.92e$-$06  &    3.95           \\
\hline
\end{tabular}
\end{table}

\begin{example}\label{ex2}
  In this test, we take the following parameters
\begin{align*}
  \upsilon_1&=0.3,\;\; \eta_1=0.5,\;\; \kappa_1=-\frac{\upsilon_1(3\sqrt{1+4\upsilon_1^2}-1)}{2(2+9\upsilon_1^2)},\;\;\zeta_1=-1,\;\; \delta_1 = \kappa_1,\;\;\beta_1=\zeta_1,\;\; \gamma_1=0, \\
  \upsilon_2&=0.3,\;\; \eta_2=0.6,\;\; \kappa_2=-\frac{\upsilon_2(3\sqrt{1+4\upsilon_2^2}-1)}{2(2+9\upsilon_2^2)},\;\;\zeta_2=-1,\;\; \delta_2 = \kappa_2,\;\;\beta_2=\zeta_2,\;\; \gamma_2=0.
\end{align*}
\end{example}

In the computation, we use our proposed second-order method, where the computational interval is  chosen as $[-15,15]$, final time is set to be $T=1$ and the initial value is taken as
  \begin{equation}
    u(x,0) = \textrm{sech}(x)e^{2ix},\quad v(x,0) = \textrm{sech}(x)e^{2ix}.
  \end{equation}
The ``exact solution" is computed on the very fine mesh $h=1/256,\tau=1/256$.

Table~\ref{tab5} and Table~\ref{tab6} list the errors and the convergence orders  for the method (\ref{diff1})-(\ref{diff3}) with $\alpha = 1.2, 1.5, 1.8, 2.0$ in the $L^\infty$-norm. Again, these results confirm the  second-order convergence both in  time and space variables. Table~\ref{tab7} and Table~\ref{tab8} list the errors and the convergence orders  for the method (\ref{fourth1})-(\ref{Richardson}) with $\alpha = 1.2, 1.5, 1.8, 2.0$ in the $L^\infty$-norm. And these results confirm the  fourth-order convergence both in time and space variables.

\section{Conclusion}\label{section6}
In this paper, we developed a linearized implicit finite difference method for the CFGLE. The method is unconditionally stable. Moreover, a rigorous analysis of the proposed difference scheme is carried out, which includes the unique solvability, the unconditional stability, and the second-order convergence in the $L^{\infty}$-norm. Numerical tests are performed to validate our theoretical findings. Our method can be easily extended to the case with spatial fourth-order accuracy. Moreover, Richardson extrapolation is used to increase the temporal accuracy to fourth order. The method proposed in this paper can also be extended to the strongly coupled nonlinear fractional Ginzburg-Landau equation, which will be our future objective.

\section*{Acknowledgement}
Dongdong He was supported by  the president's fund-research start-up fund from the Chinese University of Hong Kong, Shenzhen. Kejia Pan was supported by the National Natural Science Foundation of China (No. 41474103), the Excellent Youth Foundation of Hunan Province of China (No. 2018JJ1042) and the Innovation-Driven Project of Central South University (No. 2018CX042).

\end{document}